\documentclass[11pt]{amsart}
\usepackage{geometry} % see geometry.pdf on how to lay out the page. There's lots.
\usepackage{pgfplots}
\usepackage{psfrag}
\usepackage{amsmath}
\usepackage{calc}
\usepackage{systeme}
\usepackage{amsfonts}
\usepackage{amsthm}
\usepackage {amssymb}
\usepackage{algorithm2e}
\usepackage[french]{babel}
\usepackage[latin1]{inputenc}
\usepackage[T1]{fontenc}
\usepackage{lmodern}

\usepackage{bbold}
\usepackage{caption}
\usepackage{subcaption}
\usepackage{color,overpic,cite,pdfsync,url}
\usepackage{amsmath,amsfonts,amssymb,amsthm,bbm,bm}
\usepackage{epstopdf,float}

\usepackage[numbers]{natbib}

\newtheorem{Theorem}{Theorem}[section]

\usepackage[foot]{amsaddr}

\usepackage{appendix}
\usepackage{chngcntr}
\usepackage{etoolbox}

\title{Strategic advantages in mean field games with a major player }
\author{Charles Bertucci$^1$, Jean-Michel Lasry$^2$, Pierre-Louis Lions$^{2,3}$}
\date{}

\thanks{$^1$ : CMAP, \'Ecole Polytechnique, CNRS, 91128 Palaiseau\\ $^2$ :Universit\'e Paris-Dauphine, PSL Research University,UMR 7534, CEREMADE, 75016 Paris, France\\
$^3$ : Coll\`ege de France, 3 rue d'Ulm, 75005, Paris, France
}

\begin{document}
\maketitle

%\section*{Acknowledgment}

%The preferred spelling of the word ``acknowledgment'' in America is without an ``e'' after the ``g''. Avoid the stilted expression ``one of us (R. B. G.) thanks $\ldots$''. Instead, try ``R. B. G. thanks$\ldots$''. Put sponsor acknowledgments in the unnumbered footnote on the first page.
\begin{abstract}
%% Text of abstract
Cette note porte sur une problématique de modélisation issue de la th\'eorie des jeux \`a champ moyen. On montre comment il est possible de mod\'eliser des jeux \`a champ moyen avec un agent majoritaire qui a un avantage stratégique, tout en restant dans un cas o\`u on ne consid\`ere que des strat\'egies markoviennes en boucles fermées pour tous les joueurs. Nous illustrons ce fait autour de trois exemples.\\
\\
\noindent ABSTRACT :
This note is concerned with a modeling question arising from the mean field games theory. We show how to model mean field games involving a major player which has a strategic advantage, while only allowing closed loop markovian strategies for all the players. We illustrate this property through three examples.
\end{abstract}

\section{Introduction}
In this note we present different mean field games (MFG for short) models involving a major player who has some strategical advantages (in a sense to be detailed later on). We introduce systems of partial differential equations (PDE) which characterize Nash equilibria of MFG with a major player in the three following cases : when the decision horizon of the major player is much larger than the ones of the non-atomic players,  when the major player can control the set of actions of the non-atomic players and finally when the major player can decide to stop the game. 

MFG models have been introduced by the last two authors first in the case without common noise for the players \cite{lasry2007mean} and next in the general case with the introduction of the so-called master equation (see \citep{lions2007cours}), which is in general an infinite dimensional nonlinear partial differential equation. Let us also recall that the essential structure conditions (monotonicity) and most of the existing mathematical tools have also been introduced by the last two authors (see \citep{lasry2007mean,lions2007cours}). Let us mention that the particular case of MFG models with no common noise was independently considered in \citep{huang2006large} and that some particular cases have been discussed previously in the Economics literature (anonymous games in the discrete time case without common noise, or a heuristic description in a Macroeconomics example in \citep{krusell1998income}). Finally, since their introduction, there is now a huge scientific literature on MFG concerning the mathematical theory and also many applications to Economics, Finance, Social Sciences, Communication Networks, Engineering Sciences, Computer Sciences, etc. We refer the reader to the online courses \citep{lions2007cours}, to the book \citep{carmona2018probabilistic} (and the references therein) for further information on the mathematical theory.

In many of the applicative situations in which the MFG theory is of interest, a major player is also part of the game. Let us mention for instance the presence of regulators in economic situation or a cartel facing a fringe of small producers in some markets, such as it is the case in the oil production one. Thus, the study of MFG involving a major player is of the utmost importance, both in terms of the mathematical questions raised and in terms of applications. MFG involving a major player have been considered by several authors, let us mention for instance \citep{bensoussan2016mean,carmona2018probabilistic,lasry2018mean,cardaliaguet2018remarks}. Different modeling points of view have been adopted in the literature, some using different notions than Nash equilibria to capture the behavior that one can expect to observe in practical situations. Let us mention that \citep{cardaliaguet2018remarks} unifies in some sense the different works involving Nash equilibria. A fundamental question is to analyse to which extent the major player dominates the game. Several authors have investigated the possibility of capturing the dominance of the major player by studying other notions of equilibria than the one of Nash, such as Stackelberg equilibria. In this note, we present how we can avoid such delicate (especially in the time dependent setting) notions by giving an inherent advantage to the major player (action set, time horizon...). We work in the setting of \citep{lasry2018mean}. In particular we are interested with closed-looped Nash equilibria that we are able to characterize with systems of PDE.

\subsubsection*{A remark on terminology}
In the present note, the term "small player", used to characterize the players, is given the same meaning as the terminology "non-atomic player", introduced by R. Aumann in Game Theory. Both terms refer to players who cannot individually have an impact on the state or the payoff of the other players, as those quantities only depend on mean field terms.

\section{Underlying model and notations}
We consider a so-called major player interacting with an infinite number of indistinguishable, indentical and non-atomic agents. The time is denoted by $t$. The state of the major player is described by $y\in \mathbb{R}^d$ and his value function is denoted by $\varphi$. The state space of the non-atomic agents is $\{1;...; k\}$ and the solution of the associated master equation is $U$. The unnormalized histogram describing the repartition of non-atomic players is $x \in \mathbb{R}^k$. A typical system describing a MFG with a major player is :
\begin{equation}\label{ex}
\begin{cases}
\partial_t \varphi + F(x,y,U,\nabla_y\varphi,\alpha^*) + A(x,y,U,\alpha^*)\cdot \nabla_x \varphi - \nu \Delta_y \varphi + µ \varphi = 0;\\
\partial_t U + (A(x,y,U,\alpha^*)\cdot \nabla_x)U + \alpha^* \cdot \nabla_y U - \nu \Delta_y U + \lambda U = B(x,y,U,\alpha^*);\\
\alpha^* = \partial_p F(x,y,U,\nabla_y \varphi, \alpha^*) ;\\
U|_{t = 0} = U_0 ; \varphi |_{t = 0} = \varphi_0;
\end{cases}
\end{equation}
for $x \in \mathbb{R}^k, y \in \mathbb{R}^d, t \geq 0$, where $F$ is the Hamiltonian of the major player, $\alpha^*$ its control, $µ \geq 0$ its intertemporal preference rate and $\nu > 0$ a parameter describing the intensity of a noise associated to the major player. The terms $A$ and $B$ are modeling the MFG interactions and $\lambda$ is the intertemporal preference rate of the non-atomic players. Let us note that the time has been reversed to simplify notations.

In the game associated to (\ref{ex}), the players interact both with their state and their control, the major player perceiving of course only mean field terms. We only indicate this general framework to set some notations. Of course many different settings can be addressed, some of them shall be in the rest of this note.

\section{Different intertemporal preference rates for the major player and the crowd}
In many practical situations, the major player and the non-atomic players have different intertemporal preference rates. This can be due for instance to internal objectives or to differences in cash reserves leading to the urgency of obtaining results for a given party (a cartel of major players versus a fringe of non-cooperative small producers). Here we indicate the limit of systems such as (\ref{ex}) when the non-atomic players become more and more short-sighted, i.e. the limit $\lambda \to \infty$. In \citep{bertucci2019some}, we studied the same type of limit without a major player. In the present case, the same type of limits occurs : $U \to 0 $ as $\lambda \to \infty$ and the system reduces to 
\begin{equation}\label{myop}
\begin{cases}
\partial_t \varphi + F(x,y,0,\nabla_y\varphi,\alpha^*) + A(x,y,0,\alpha^*)\cdot \nabla_x \varphi - \nu \Delta_y \varphi + µ \varphi = 0;\\
\alpha^* = \partial_p F(x,y,0,\nabla_y \varphi, \alpha^*) ;\\
\varphi |_{t = 0} = \varphi_0;
\end{cases}
\end{equation}
for $x \in \mathbb{R}^k, y \in \mathbb{R}^d, t \geq 0$. The following holds :
\begin{Theorem}
Under the assumptions of theorem 3.2 of \citep{lasry2018mean}, there exists a time $T_0 > 0$ such that, for any $0 < \epsilon < T_0$, extracting a subsequence if necessary, the sequence of solutions $(\varphi_{\lambda},U_{\lambda})_{\lambda > 0}$ of (\ref{ex}) converges locally uniformly in $(\epsilon,T_0)\times \mathbb{R}^{k}\times \mathbb{R}^d$ toward $(\varphi, 0)$ as $\lambda$ goes to infinity, where $\varphi$ is a solution of (\ref{myop}).
\end{Theorem}

In the limit the small players do not anticipate any longer. Their behavior can now be described by what is referred in the literature as agent-based models or population games, but they are now interacting with a rational player with anticipations. Let us note that in the limit, the controls of the small players are still optimal, but for a different problem, which has zero duration. We also believe that this type of model is well-suited to study Stackelberg like equilibria. Indeed, as the time dependence disappears for the small players, such concepts become easier to manipulate. For instance we can think that this model may be used to study games in which the time scale of the non-atomic players is the day, and the time scale of the major player involves one year or more.

Finally, let us remark that such models seem particularly adapted to study questions related to climate change and to limiting its causes. Indeed a (or several) regulator/ major player should have a significantly longer time scale in mind, compared to the one of the crowd of the different agents who are impacting the environment.

\section{The major player can take the "lead"}
In this section we are interested in games in which the actions of the non-atomic players can be directly impacted by the state of the game or the actions of the major player. In particular, we emphasize the fact that, while continuing to work in the setting of closed-loop markovian Nash equilibria, we can study games in which the controls of the small players can be inhibited by the major player. This can be obtained without modifying the key structure of the system (\ref{ex}).

To set ideas on an example, take the case of a continuous state space for the non-atomic players, when the system takes the form :
\begin{equation}
\begin{cases}
\partial_t \varphi + F(m,y,U[m],\nabla_y\varphi,\alpha^*) + \langle \frac{\partial \varphi}{\partial m},- \text{div}(\frac{\partial H}{\partial p}m)\rangle - \nu \Delta_y \varphi + µ \varphi = 0;\\
\partial_t U + H(x,y,m,\nabla_x U, \alpha^*)+  \langle \frac{\partial U}{\partial m},- \text{div}(\frac{\partial H}{\partial p}m)\rangle + \alpha^* \cdot \nabla_y U - \nu \Delta_y U + \lambda U = 0;\\
\alpha^* = \partial_p F(m,y,0,\nabla_y \varphi, \alpha^*) ;\\
U|_{t = 0} = U_0 ; \varphi |_{t = 0} = \varphi_0.
\end{cases}
\end{equation}
We refer to \citep{lasry2018mean} for the interpretation of the different terms in the previous system. Let us assume that the Hamiltonian of the non-atomic players $H$ is such that there exists $\tilde{H}$ and $V$ such that
\begin{equation}
H(x,y,m,p, \alpha) = \alpha\tilde{H}(x,y,m,p) + V \cdot p;
\end{equation}
then the major player can choose $\alpha = 0$ and force the small players to move along the vector field $V$. Another interesting situation is the case in which the Hamiltonian $H$ has the form 
\begin{equation}
H(x,y,m,p, \alpha) = G(x,y,m)\tilde{H}(x,y,m,p) + \alpha \cdot p.
\end{equation}
In such a situation, on the set $\{G = 0\}$, the small players cannot control their position and they move only according to $\alpha$.

Applications involving such modeling aspects are numerous. Let us mention the important case of a regulator which has the possibility to impact the effect of the controls of the small players or which can take control, under certain circumstances, of the evolution of the state of the small players.

\section{Optimal stopping by the major player}
In this section, we consider games in which the major player can decide to stop the game (for everyone). More generally, we could consider impulse or switching controls for the major player. For instance, we could consider games in which the major player can switch between two states which are both associated to a system of the form of (\ref{ex}). As the concepts involved in the following analysis are the same, we restrict ourselves to the case of optimal stopping to simplify notations.

This type of setting is of course reminiscent of a referee whistling the end of a game and applications of this type of models include for instance the closure of a market by a regulator.

We denote by $\psi(x,y,m)$ the stopping cost of the major player and by $\bar{U}$ the cost the non-atomic players face once the major player has stopped the game. The action set of the major player being non-convex, we cannot expect to obtain existence of Nash equilibria without considering mixed strategies. Following \citep{bertucci2018optimal,bertucci2018fokker}, we start by writing a penalized system. At this penalized level, we assume that the major player cannot decide to leave the game instantaneously, but control the intensity of an inhomogeneous Poisson process on some probabilistic space $(\Omega, \mathcal{A},\mathbb{P})$, and that this intensity has to be bounded by $\epsilon^{-1}$ for some $\epsilon > 0$. The system which describes such a game is
\begin{equation}\label{penal}
\begin{cases}
\partial_t \varphi + F(x,y,U,\nabla_y\varphi,\alpha^*) + A(x,y,U,\alpha^*)\cdot \nabla_x \varphi - \nu \Delta_y \varphi + µ \varphi + \beta^*(\varphi - \psi)^+= 0;\\
\partial_t U + (A(x,y,U,\alpha^*)\cdot \nabla_x)U + \alpha^* \cdot \nabla_y U - \nu \Delta_y U + \lambda U + \beta^*(U - \bar{U})= B(x,y,U,\alpha^*);\\
\alpha^* = \partial_p F(x,y,U,\nabla_y \varphi, \alpha^*) ;\\
0\leq \beta^*\leq  \epsilon^{-1} ; \beta^* \equiv \epsilon^{-1} \text{ on } \{\varphi > \psi\} ; \beta^* = 0 \text{ on } \{\varphi < \psi\};\\
U|_{t = 0} = U_0 ; \varphi |_{t = 0} = \varphi_0;
\end{cases}
\end{equation}
for $x \in \mathbb{R}^k, y \in \mathbb{R}^d, t \geq 0$. In this system, $\beta^*$ is the optimal intensity of the Poisson process controlled by the major player. In the first equation, the term $\beta^*(\varphi - \psi)^+$ is the standard penalization term for optimal stopping. In the second equation, the term $\beta^*(U - \bar{U})$ stands for the fact that the small players are anticipating that at the rate $\beta^*$, the game may stop, in which case they will face the cost $\bar{U}$. Existence results for this system can be obtained using arguments introduced in \citep{lasry2018mean} and fixed point argument on correspondances.
\begin{Theorem}
Under the assumption of theorem $3.2$ of \citep{lasry2018mean} together with the assumption that $\psi$ and $\bar{U}$ are smooth functions, for any $\epsilon > 0$, there exists $T_0> 0$ such that there exists a solution $(\varphi, U, \alpha^*, \beta^*)$ of (\ref{penal}) for $T< T_0$.
\end{Theorem}
The passage to the limit $\epsilon \to 0$ in (\ref{penal}) yields the following system 
\begin{equation}\label{stop}
\begin{cases}
\max\left\{\varphi - \psi ; \partial_t \varphi + F(x,y,U,\nabla_y\varphi,\alpha^*) + A(x,y,U,\alpha^*)\cdot \nabla_x \varphi - \nu \Delta_y \varphi + µ \varphi \right\} = 0;\\
\begin{aligned}\partial_t U &+ (A(x,y,U,\alpha^*)\cdot \nabla_x)U + \alpha^* \cdot \nabla_y U - \nu \Delta_y U\\
 &+ \lambda U + \beta^*(U - \bar{U})= B(x,y,U,\alpha^*) \text{ on } \{ \beta^* < + \infty\};\end{aligned}\\
0\leq \beta^*\leq  + \infty ;\\
U = \bar{U} \text{ on } \{\partial_t \varphi + F(U,\nabla_y\varphi,\alpha^*) + A\cdot \nabla_x \varphi - \nu \Delta_y \varphi + µ \varphi < 0\} ;\\
\beta^* = 0 \text{ on } \{\varphi < \psi\};\\
U|_{t = 0} = U_0 ; \varphi |_{t = 0} = \varphi_0;
\end{cases}
\end{equation}
Let us make the following remarks : (i) The fact that there exists a non-empty set $\{ 0 < \beta^* < + \infty\}$ (which makes this formulation quite complex) is in general unavoidable. It arises from the fact that we are forced to consider Nash equilibria in mixed strategies for the major player to obtain existence. (ii) Under some structural assumptions on $F$, $A$ and $\psi$, it can be shown that Nash equilibria in pure strategies exist. In such cases we can characterize $(\varphi, U)$ simply with the fact that $\varphi$ solves the obstacle problem, $U$ the master equation on $\{\varphi < \psi\}$ and $U = \bar{U}$ on $\{\varphi = \psi\}$. (iii) For (\ref{stop}) to have classical solutions, some compatibility assumptions have to be made on $\bar{U}$, namely that it is a solution of the master equation.

\section*{Acknowledgments}
The second and third authors have been partially supported by the Chair FDD (Institut Louis Bachelier). The third author has been partially supported by the Air Force Office for Scientific Research grant FA9550-18-1-0494 and the Office for Naval Research grant N000141712095.

\bibliographystyle{plainnat}
\bibliography{bib1}
\end{document}